\documentclass[12pt,twoside]{article}
\usepackage{amsxtra,amssymb,amsthm,amsmath,latexsym}

\theoremstyle{plain}

\def\oH{\buildrel\circ\over H}
\def\oH1{\buildrel\circ\over H\kern-.02in{}^1}
\def\oH1{\buildrel\circ\over H\kern-.02in{}^1}

\def\d{\delta}
\def\ep{\epsilon}
\def\la{\lambda}

\begin{document}

\title{
 Discrepancy principle for the dynamical systems method\\
\footnotesize{ Key words: ill-posed
problems, dynamical systems method (DSM), discrepancy principle,
evolution equations.}\\
\footnotesize{Math subject classification: 34R30,  35R25, 35R30, 37C35, 37L05,
37N30, 47A52, 47J06, 65M30, 65N21; }\\ 
\footnotesize{PACS 02.30.-f, 02.30.Tb,
02.30.Zz,02.60Lj, 02.60.Nm, 02.70.Pt, 05.45.-a}\\}

\author{
A.G. Ramm\\
  Mathematics Department, Kansas State University, \\
 Manhattan, KS 66506-2602, USA\\
ramm@math.ksu.edu
}

\date{}

\maketitle

\centerline{ Abstract} Assume that
$$
Au=f,\quad (1)
$$
is a solvable linear equation in a Hilbert space, $||A||<\infty$,
and $R(A)$ is not closed, so problem (1) is ill-posed. Here $R(A)$ is the
range of the linear operator $A$. A DSM (dynamical systems method)
for solving (1), consists of solving the following Cauchy problem:
$$
\dot u= -u +(B+\ep(t))^{-1}A^*f, \quad u(0)=u_0, \quad (2)
$$
where $B:=A^*A$, $\dot u:=\frac {du}{dt}$, $u_0$ is arbitrary, and
$\ep(t)>0$ is a continuously differentiable function, monotonically
decaying to zero as $t\to \infty$.
A.G.Ramm has proved that, for any $u_0$, problem (2) has a
unique solution for all $t>0$, there exists $y:=w(\infty):=\lim_{t\to
\infty}u(t)$, $Ay=f$, and $y$ is the unique minimal-norm solution to
(1). If $f_\d$ is given, such that $||f-f_\d||\leq \d$, then
$u_\d(t)$ is defined as the solution to (2) with $f$ replaced by $f_\d$.
The stopping time is defined as a number $t_\d$
such that $\lim_{\d \to 0}||u_\d (t_\d)-y||=0$, and $\lim_{\d \to
0}t_\d=\infty$. A discrepancy principle is proposed and proved in this
paper. This principle yields $t_\d$ as the unique solution to the
equation:
$$
||A(B+\ep(t))^{-1}A^*f_\d -f_\d||=\d, \quad  (3)
$$
where it is assumed that $||f_\d||>\d$ and $f_\d\perp N(A^*)$.
For nonlinear monotone $A$ a discrepancy principle is formulated and 
justified.

\section{Introduction and statement of the result. }

Assume that
$$
Au=f,
 \eqno{(1.1)}
$$
is a solvable linear equation in a Hilbert space, $||A||<\infty$,
and $R(A)$ is not closed, so problem (1.1) is ill-posed. Here $R(A)$ is 
the
range of the linear operator $A$. Without loss of generality, assume that
$||A||\leq 1$.
Let $y$ be the unique minimal-norm solution to (1.1).
A solvable equation (1.1) is equivalent to
$$
Bu=A^*f, \quad B:=A^*A,
 \eqno{(1.2)}
$$
where $A^*$ is the operator adjoint to $A$. One has
$N(B)=N(A):=N:=\{v: Av=0\}$. Let $Q:=AA^*$, and $a>0$ be a number. Then
$||B||\leq 1$, $||Q||\leq 1$, and
$(B+a)^{-1}A^*=A^*(Q+a)^{-1}$, as one easily checks. Denote by $E_\la$
and $F_\la$ the resolutions of the identity of $B$ and $Q$, respectively.

Let $\ep(t)$ be a monotone, decreasing function, 
$$
\ep(t)>0,\quad \lim_{t\to 0} \ep(t)=0,\quad  \lim_{t\to \infty} \sup 
_{\frac t 2\leq s\leq t} |\dot \ep(s)|\ep^{-2}(t)=0. 
$$
A DSM (dynamical systems method)
for solving (1.1), consists of solving the following Cauchy problem:
$$
\dot u= -u +(B+\ep(t))^{-1}A^*f, \quad u(0)=u_0,\quad \dot u:=\frac
{du}{dt},
 \eqno{(1.3)}
$$
where $u_0$ is arbitrary, and
proving that, for any $u_0$, problem (1.3)
has a unique solution for all $t>0$, there exists
$y:=u(\infty):=\lim_{t\to \infty}u(t)$, $Ay=f$, and $y$ is the unique
minimal-norm solution to (1.1). These results are proved in [3].
If $f_\d$ is given, such that
$||f-f_\d||\leq \d$, then
$u_\d(t)$ is defined as the solution to (1.3) with $f$ replaced by $f_\d$.
The stopping time is defined as a number $t_\d$
such that $\lim_{\d \to 0}||u_\d (t_\d)-y||=0$, and $\lim_{\d \to
0}t_\d=\infty$. A discrepancy principle for choosing $t_\d$ is proposed
and proved in this paper.

Let us assume $f_\d \perp N(A^*)$. Then this discrepancy principle,
yields $t_\d$ as the unique solution to the equation:
$$
||A(B+\ep(t))^{-1}A^*f_\d -f_\d||=\d, \quad ||f_\d||>\d,
 \eqno{(1.4)}
$$
and we prove that
$$
\lim_{\d \to 0}||u_\d (t_\d)-y||=0, \quad lim_{\d \to 0}t_\d=\infty.
 \eqno{(1.5)}
$$

 Thus, our basic result is:

{\bf Theorem 1.1.} {\it If $A$ is a bounded linear operator in a Hilbert
space $H$, equation (1.1) is solvable, $||f_\d||>\d$,
$f_\d \perp N(A^*)$, and $\ep(t)$
satisfies the assumptions stated above, then equation (1.4) has a
unique solution $t_\d$, and (1.5) holds, where $y$ is the unique 
minimal-norm solution to (1.1).}

In Section 2 a proof of this theorem is given, and a discrepancy principle
is proved for equation (1.1) with monotone, nonlinear, continuous
operator $A$.
For variational regularization the discrepancy principle was proposed by
Morozov [2], see also [1] and [4]. In [5] a general regularization method
is proposed for a wide class of nonlinear ill-posed problems.

\section{Proofs. }

First, let us prove two lemmas.
By $\rightharpoonup$ and $\rightarrow$ we denote the weak and strong
convergence in $H$, respectively.

{\bf Lemma 1.} {\it If $w_n\rightharpoonup y$, and $\limsup_{n\to
\infty}||w_n||\leq ||y||$, then $w_n\rightarrow y$.}

{\bf Lemma 2.} {\it If $\ep(t)>0$, $\dot \ep<0$, and $\lim_{t\to \infty}
\dot \ep(t) \ep^{-2}(t)=0$, then $\lim_{t\to \infty}e^{-t}\ep^{-1}(t)=0$.}

{\bf Proof of Lemma 1.} If $w_n\rightharpoonup y$, then $||y||\leq
\liminf_{n\to \infty} ||w_n||$. This, and the inequality $\limsup_{n\to
\infty}||w_n||\leq ||y||$ yield $\lim_{n\to \infty} ||w_n||=||y||$.
Thus, $\lim_{n\to \infty} ||w_n-y||^2=\lim_{n\to \infty} [||w_n||^2
+||y||^2-2\Re (w_n,y)]=0$. $\Box$

{\bf Proof of Lemma 2.} Our assumptions imply $\frac{d \ep^{-1}(t)}{dt}
\leq c$, where $c=const>0$. Thus $\ep^{-1}(t)\leq ct +c_0$, where $c_0>0$ is a
constant. The conclusion of Lemma 2 follows. $\Box$

{\bf Remark 1.} An example of $\ep(t)$, satisfying all the assumptions
imposed in the Introduction, is $\ep(t)=c_1(c_0+t)^{-b}$, where
$c_0,c_1>0 $ are positive constants, and $b\in (0,1)$ is a constant.

The proof of Theorem 1.1 consists of two steps.

Step 1: prove that (1.4) has a unique solution $t_\d$, and $\lim_{\d \to 
0}
t_\d=\infty$.

Step 2: prove (1.5).

Step 1. Write $\ep$ for $\ep(t)$ and rewrite (1.4) as:
$$
\d^2= ||A(B+\ep)^{-1}A^*f_\d -f_\d||^2=||[Q(Q+\ep)^{-1}-I]f_\d||^2=
$$
$$=\ep^2 \int_0^1(\ep+\la)^{-2}d\rho(\la):=h(\d,\ep).
\eqno{(2.1)}
$$
Here $Q=AA^*$, $\rho=(F_\la f_\d, f_\d)$, $||Q||\leq 1$.
One has $h(\infty, \d)=\int_0^1d\rho=||f_\d||^2>\d^2$,
and $\lim_{\ep \to 0} h(\ep, \d)=0,$ by the dominant convergence theorem,
provided that $\lim_{\gamma \to 0}\int_0^\gamma d\rho=0$.
The last relation holds if
$f_\d \perp N(A^*)=N(Q)$.  This assumption is natural, because
$f_\d$ enters under the sign of $A^*$ in the definition of
$w_\d(s)$ in the argument given in Step 2.
Thus, $h(\d, \infty)>\d^2$ and $h(\d,0)=0<\d^2$. By the continuity of
$h$ as a function of $\ep$ one concludes that there exists a solution
$\ep:=g(\d)>0$ to the equation $h(\d, \ep)=\d^2$.
Since $h$ is a monotone increasing function of $\ep$, this solution 
is
unique. Because $\ep(t)$ is monotone decreasing, the equation
$\ep(t)=g(\d)$ defines a unique root $t_\d$, and
$\lim_{\d \to 0}t_\d=\infty$. Step 1 is done.

Step 2. To prove (1.5), write
$u_\d(t)=u_0e^{-t}+\int_0^t e^{-(t-s)}w_\d(s)ds$, where here and below
$t:=t_\d$ is defined in Step 1, and
$w_\d(s):=(B+\ep(s))^{-1}A^*f_\d$. Because $f_\d$ enters under the sign
of $A^*$,  one may assume that $f_\d \perp N(A^*)= N(Q)$.
This assumption has been used in Step 1.

We have $u_\d=j_1 +j$, where $j_1:=u_0e^{-t}+\int_0^{t/2}
e^{-(t-s)}w_\d(s)ds$, and
$j:=\int_{t/2}^t e^{-(t-s)}w_\d(s)ds.$
One has $||w_\d(s)|\leq ||f_\d||\ep^{-1}(s)$, because $||A^*||\leq 1,$
and $||(B+\ep)^{-1}||\leq \ep^{-1}$. By Lemma 2, 
$$||j_1||\leq
||u_0||e^{-t} +e^{-t/2}\ep^{-1}(t/2)||f_\d||\to 0\quad as\quad t\to 
\infty.
$$
Furthermore, 
$$||j-y||\leq ||w_\d(t) \int_{t/2}^te^{-(t-s)} ds-y||+
||\int_{t/2}^t e^{-(t-s)}[w_\d(s)-w_\d(t)] ds ||:=J_1+J_2.
$$
One has: 
$$J_1=||w_\d(t)-y-w_\d(t)e^{-t/2}||\to 0\quad as \quad t\to \infty,$$
because, as we prove below,
$$
\lim_{\d \to 0} ||w_\d(t_\d)-y||=0.
\eqno{(2.2)}
$$
To prove that $J_2\to 0$ as $t\to \infty$, we estimate
$$||w_\d(s)-w_\d(t)||=||(B+\ep(s))^{-1}(\ep(s)-\ep(t))(B+\ep(t))^{-1}A^*f_\d||
\leq ||f_\d||\ep^{-2}(t)|\dot \ep (\xi)|(t-s),$$
where $t/2\leq s\leq t$, and $\xi$ is an intermediate
point in the Lagrange formula. Thus,
$$J_2\leq ||f_\d|| \int^t_{t/2} e^{-(t-s)}(t-s)ds \ep^{-2}(t)
\sup_{\frac t 2 \leq \xi \leq t}|\dot \ep (\xi)|\to 0 \quad  as \quad 
t\to \infty.$$

Let us prove (2.2). The element $w_\d:=w_\d(t_\d)$ is
the minimizer of the problem: $||Aw-f_\d||^2+\ep ||w||^2=\min$,
$\ep:=\ep(t),\,\, t:=t_\d.$ Thus,
$$||Aw_\d-f_\d||^2+\ep ||w_\d||^2\leq ||A(y)-f_\d||^2+\ep ||y||^2= \d^2 
+\ep||y||^2.
$$
 So, $||w_\d||\leq ||y||$, because of (1.4). Therefore, one may assume
that $w_\d \rightharpoonup W$  as $\d \to 0$, and, as we prove below,
$W=y$.
Now, Lemma 1 implies that (2.2) holds.

Let us prove that $W=y$. It follows from (1.4) that $Bw_\d\to A^*f$
 as $\d \to 0$. This
and $w_\d\rightharpoonup W$, imply $BW=A^*f$, since $B$ is monotone and
therefore $w-$closed. We prove $w-$closedness of $B$ 
below.
Since the minimal-norm solution to (1.2) is unique, and since $||W||\leq
||y||$, it follows that $W=y$, as claimed.

Finally, we prove $w-$closedness of a monotone, hemicontinuous operator 
$B$
defined on all of $H$. We call an operator $B$ 
$w-$closed if
$w_j\rightharpoonup y$ and $Bw_j\to f$ imply $By=f$.
An operator $B$ is called monotone, if
$(B(u)-B(v),u-v)\geq 0$ for all $u,v\in D(B)$, and
hemicontinuous, if the function $t\mapsto (B(u+tv),z)$
is continuous for any $u,v,$ and $z\in H$, as a function of $t\in [0,1)$.
To prove $w-$closedness of $B$, we start with the relation:
$(B(u_j)-B(y-tz),u_j-y+tz)\geq 0$.
Let $j\to \infty$. Then $(f-B(y-tz),tz)\geq 0$, or, since $t>0$,
 $(f-B(y-tz),z)\geq 0$. Let $t\to 0$. Then, by hemicontinuity of $B$,
one gets  $(f-B(y),z)\geq 0$. Since $z$ is arbitrary, it follows that
$By=f$,
as claimed. Note, that if $B$ is continuous, it is hemicontinuous. Step 2
is done. Theorem 1.1 is proved. $\Box$

{\bf Remark 2.} Suppose that the assumption $f_\d \perp N(A^*)$ does not
hold. Let $f_\d=g+h$, where $g\perp h$, $g\in N(A^*)$. Then one can use
the discrepancy principle of the form:
$$
||A(B+\ep(t))^{-1}A^*f_\d -f_\d||=C\d, \quad ||f_\d||>C\d,\quad
C=const>1.
 \eqno{(1.4')}
$$
 In Step 1 we prove now that $h(\d, 0)< C^2\d^2,$  $h(\d, \infty)>
C^2\d^2,$ and the arguments
are similar to the given in the proof of Theorem 1.1.

{\bf Remark 3.} In this remark we prove a {\it discrepancy
principle for nonlinear, monotone, continuous operators,
defined on all of $H$.}

{\bf Theorem 2.1.} {\it Assume: i) $A$ is a monotone, continuous, 
defined on all of $H$, operator, ii) equation $A(u)=f$ is solvable,
$y$ is its minimal-norm solution, and iii) $||f_\d -f||\leq \d$,
$||A(0)-f_\d||>C\d$, where $C>1$ is a constant which can be 
chosen arbitrarily close to $1$. Define $\ep(\d)>0$ to be 
the minimal solution to the equation 
$$ ||A(u_{\d,\ep})-f_\d||=C\d,
 \eqno{(2.3)}
$$
where $u_{\d,\ep}$ is any element satisfying inequality 
$$F(u_{\d,\ep})\leq 
m+(C^2-1)\d^2,  \quad F:=||A(u)-f_\d||^2+\ep ||u||^2,
$$
 and
$m=m(\d,\ep):=inf_{u} F(u)$. Then equation (2.3) for $\ep$ has a solution
$\ep(\d)>0$. If $\ep(\d)$ is its minimal solution, and 
$u_\d:=u_{\d,\ep(\d)}$, then $\lim_{\d \to 0}||u_\d-y||=0$.}

{\bf Proof. } If $A$ is monotone, continuous, defined on all of 
$H$, 
then the set $N_f:=\{u: A(u)=f\}$ is convex and closed, so it has 
a 
unique minimal-norm element $y$.
To prove the existence of a solution to (2.3), we prove that
the function $h(\d, \ep):=||A(u_{\d,\ep})-f_\d||$ is greater than $C\d$ 
for sufficiently large $\ep$, and smaller than $C_1\d$ for sufficiently 
small $\ep$ and $C_1>C$ arbitrarily close to $C$. Because of the 
continuity of $h(\d, \ep)$ with respect to $\ep$ on $(0,\infty)$, equation 
$h(\d, \ep)=C\d$ has a solution if $C>1$ and $C<||A(0)-f_\d||$.

Let us prove the above estimates. One has
$$F(u_{\d,\ep})\leq m+(C^2-1)\d^2 \leq F(0)+(C^2-1)\d^2,$$
and 
$$F(u_{\d,\ep})\leq m+(C^2-1)\d^2 \leq F(y)+(C^2-1)\d^2= 
\ep ||y||^2+C^2\d^2.$$ 
Therefore, as $\ep \to \infty$, one gets
$||u_{\d,\ep}||\leq \frac c {\sqrt {\ep}}\to 0,$ where $c>0$ is a constant
depending on $\d$. Thus, by the continuity of $A$, one obtains
$\lim_{\ep\to \infty}h(\d, \ep)= ||A(0)-f_\d||>C\d$. 

Now, let $\ep \to 0$. Then
$h^2(\d,\ep)< \ep ||y||^2+C^2\d^2$. Thus $\liminf_{\ep \to 0}h(\d,\ep)\leq 
C\d$. Therefore equation $h(\d,\ep)=C_1\d$ has a solution $\ep(\d)>0$
if $C<C_1<||A(0)-f_\d||$, which is what we want, because $C_1>1$ can be 
taken arbitrarily close to $1$ if $C>1$ 
can be taken arbitrarily close to $1$.

Let us now prove that if $u_\d:=u_{\d,\ep(\d)}$, then $\lim_{\d\to 
0}||u_\d-y||=0$. From the estimate 
$$||A(u_\d)-f_\d||^2+\ep||u_\d||^2\leq C^2\d^2+\ep||y||^2,
$$ 
and from the equation (2.3), it follows that 
$||u_\d||\leq ||y||$. Thus, one may assume that $u_\d\rightharpoonup U$,
and from (2.3) it follows that $A(u_\d)\to f$ as $\d\to 0$. By 
$w-$closedness 
of monotone continuous operators, one gets $A(U)=f$, and from
$||u_\d||\leq ||y||$ it follows that $||U||\leq ||y||$. Because the 
minimal 
norm solution to the equation $A(u)=f$ is unique, one gets $U=y$. Thus,
$u_\d\rightharpoonup y$, and $||u_\d||\leq ||y||$. By Lemma 1, it follows 
that $\lim_{\d\to 0}||u_\d-y||=0$. Theorem 2.1 is proved. $\Box$.

\end{document}